\documentclass[12pt,reqno]{amsart}

\usepackage{fullpage,amsfonts,amsmath,amssymb,amscd,}

\input xypic

\theoremstyle{plain}
\newtheorem*{lem}{Lemma}
\newtheorem*{prop}{Proposition}

\newtheorem*{thrm}{Theorem}
\newtheorem*{cor1}{Corollary 1}
\newtheorem*{cor2}{Corollary 2}
\newtheorem*{cor3}{Corollary 3}
\newtheorem*{cor4}{Corollary 4}
\theoremstyle{definition}
\newtheorem*{examples}{Examples}
\newtheorem*{defn}{Definition}

\newtheorem*{rem}{Remark}

\newtheorem*{claim}{Claim}

\newcommand{\Spec}{\mathrm{Spec}}

\newcommand{\N}{\mathbb{N}}
\newcommand{\Z}{\mathbb{Z}}
\newcommand{\C}{\mathbb{C}}

\renewcommand{\leq}{\leqslant}
\renewcommand{\geq}{\geqslant}

\newcommand{\grA}{\mathrm{gr} ^{\mathcal{F}} A}
\newcommand{\grI}{\mathrm{gr} ^\mathcal{F} I}
\newcommand{\s}{\mathcal{S}}
\newcommand{\os}{\overline{\mathcal{S}}}
\newcommand{\var}{\mathcal{V}}
\newcommand{\grM}{\mathrm{gr}^{\mathcal{F}} M}
\newcommand{\grGM}{\mathrm{gr}^{\Gamma} M}
\newcommand{\grf}{\mathrm{gr}^{\mathcal{F}}}
\newcommand{\F}{\mathcal{F}}
\newcommand{\fv}{\mathbf{v}}
\newcommand{\Ug}{\mathcal{U}(\mathfrak{g})}

\title{The associated variety of a Poisson prime ideal}
\author{Maurizio Martino}

\email{mma@maths.gla.ac.uk}

\address{Department of Mathematics, University
of Glasgow, Glasgow, G12 8QW, U.K.}

\date{11th May 2004}

\begin{document}

\begin{abstract}
We prove that the associated variety of a Poisson prime ideal of
the centre of a symplectic reflection algebra at parameter $t = 0$
is irreducible.
\end{abstract}

\maketitle

\section{Introduction}

\subsection{} The study of the primitive ideals of
an algebra is an approximation of its representation theory. One
model case where this study has been fruitful is that of the
enveloping algebra of a complex semisimple Lie algebra. The
history of our result can be traced back to the well known theorem
by Borho and Brylinski \cite{BB} and by Joseph \cite{J} stating
that in the above context, the associated variety of a primitive
ideal is irreducible, and in fact that result can be read off as a
corollary of our main theorem.

\subsection{} The symplectic reflection algebras of Etingof and
Ginzburg \cite{EG} are an interesting class of algebras with
applications to integrable systems, invariant theory and geometry.
Their behaviour varies according to a parameter $t$. When $t \neq
0$ they have trivial centres, while when $t = 0$ they have large
centres and the geometry of the centres plays a leading role.

Let $\fv$ be an even dimensional complex vector space with
symplectic form $\omega$, and let $G$ be a finite subgroup of the
symplectic group of $\fv$. Let the triple $(\fv, \omega, G)$ be
indecomposable. The symplectic reflection algebras, $H_{t,
\bf{c}}$, are deformations of the skew group ring, $\C[\fv]
* G$, and their spherical subalgebras, $eH_{t, \bf{c}}e$, are
deformations of $\C[\fv]^G$, the ring of $G$-invariants.

\subsection{} It was proved by Ginzburg \cite[Theorem 2.1]{Gi} that
when $t \neq 0$ the associated variety of a primitive ideal of
$eH_{t, \bf{c}}e$ is irreducible. His method was to generalise the
Lie theoretic result using the ideas of Poisson geometry and
symplectic leaves. We extend this result to the case when $t = 0$.
Here the algebra $eH_{0, \bf{c}}e$ is commutative and is
isomorphic to the centre, $Z_{0, \bf{c}}$, of $H_{0, \bf{c}}$. In
fact, the centre has the structure of a Poisson algebra. It was
shown by Brown and Gordon in \cite{BrGo} that it is the Poisson
prime ideals of $Z_{0, \bf{c}}$ which provide the natural first
step in understanding the finite dimensional representation theory
of $H_{0, \bf{c}}$. Our main result is the following.\medskip
\begin{center}
\textit{The associated variety of a Poisson prime ideal of $Z_{0,
\bf{c}}$ is irreducible.}\\
\end{center}

\footnotetext[1]{ The research described here will form part of
the author's PhD thesis at the University of Glasgow. Part of the
work for the paper was done while the author participated in a
workshop at the University of Washington in August 2003, supported
by Leverhulme Research Interchange Grant F/00158/X. The author
thanks Ken Brown and Iain Gordon for suggesting this problem and
for their advice and encouragement.}

\medskip

In fact, the result we prove (Theorem \ref{Ithrm}) is somewhat
more general than this, and includes Ginzburg's result as a
special case. Our proof is modelled on Ginzburg's proof of
\cite[Theorem 2.1]{Gi}, which in turn is based on a proof by Vogan
\cite[$\S 3-4$]{Vo} of the result for enveloping algebras.

\subsection{} As discussed in \ref{rep}, our hope is that this will
allow some kind of description of irreducible finite dimensional
representations of $H_{0, \bf{c}}$ by subgroups of $G$. We hope to
study this further in later work.

Our paper is organised as follows. In $\S \ref{main}$ we introduce
basic definitions and state the main theorem; we discuss
applications in $\S \ref{appl}$. We prove the main theorem in the
remaining sections: $\S \ref{Prelim}$ states a number of
preliminary results which we use in the proof, which is given in
$\S \ref{proof}$.

\section{The Main Theorem}\label{main}

\subsection{Poisson structures}\label{Poisson}

\begin{defn}
Let $R$ be an affine commutative $\C$-algebra. We say that $R$ is
a Poisson algebra if there exists a non-trivial Poisson bracket
$\{-,-\}: R \times R \to R$. That is, $\{-,-\}$ is a non-zero
skew-symmetric bilinear map, and for all $x,y,z \in R$, $\{xy,z\}
= x\{y,z\} + \{x,z\}y$ and $ \{x,\{y,z\}\} = \{\{x,y\},z\} +
\{y,\{x,z\}\}.$ We shall say that an affine variety over $\C$ is
\textit{Poisson} if its coordinate ring is a Poisson algebra.
\end{defn}

Let $(R, \{-,-\}_R)$ and $(S, \{-,-\}_S)$ be Poisson algebras. We
say that a map $ \psi : R \to S$ is a $\it Poisson\ homomorphism$
if it is an algebra homomorphism such that for all $x, y \in R$,
$\psi (\{x,y\}_R) = \{\psi(x), \psi(y)\}_S$. When $R =
\bigoplus_{i=0}^{\infty} R_i$ is a graded Poisson algebra we shall
say that $\{-,-\}$ \textit{has degree l} if for
all $i$ and $j$, $\{R_i , R_j\} \subseteq R_{i+j+l}$, and $l$ is
the minimal integer for which this is true.\\

Let $R$ be a Poisson algebra and fix an algebra generating set
$\{a_1, \dots a_k\}.$ For a closed point $\mathbf{m} \in
\rm{Spec}\ R$ we define the rank of the Poisson structure at
$\mathbf{m}$ to be the rank of the matrix $( \{a_i, a_j\} +
\mathbf{m}) \in M_{k}(\C)$. It is independent of the choice of
generators.

\begin{defn}
The \textit{symplectic leaf} $\mathcal{S}(\mathbf{m})$ containing
a closed point $\mathbf{m}$ of $\rm{Spec}\ R$ is the maximal
connected complex analytic manifold in $\rm{Spec}\ R$ such that
$\mathbf{m} \in \rm{Spec}\ R$ and the rank of each closed point in
$\mathcal{S}(\mathbf{m})$ equals the dimension of
$\mathcal{S}(\mathbf{m}).$
\end{defn}

The symplectic leaves of $\rm{Spec}\ R$ are related to certain
ideals of $R$.

\begin{defn} Let $I$ be an ideal of a Poisson algebra $R$.
Then $I$ is a \textit{Poisson ideal} if $\{R,I\} \subseteq I$.
\end{defn}

It was shown in \cite[Proposition 1.3]{W} that when $\rm{Spec}\ R$
is smooth there exists a stratification of $\rm{Spec}\ R$ by
symplectic leaves. One can extend this as in \cite[$\S 3.5$]{BrGo}
to show that for any Poisson algebra, $R$, there exists a
stratification of $\rm{Spec}\ R$ by symplectic leaves. For any
symplectic leaf $\mathcal{S}$ ($= \mathcal{S}(\mathbf{m})$ for
some closed point $\mathbf{m} \in \rm{Spec}\ R$) the closure,
$\overline{\mathcal{S}}$, of $\mathcal{S}$ in $\rm{Spec}\ R$ is
the closed subset, $\mathcal{V}(P)$, of $\rm{Spec}\ R$ whose
defining ideal, $P$, is a Poisson prime ideal of $R$ (see
\cite[Lemma 3.5]{BrGo}).

If $\rm{Spec}\ R$ is a finite union of symplectic leaves, we say
that $\rm{Spec}\ R$ \textit{has finitely many symplectic leaves}
(or $\rm{Spec}\ R$ has FMSL). One consequence of this is that then
the notion of a symplectic leaf becomes algebraic (as opposed to
analytic) in the sense that each leaf is a locally closed subset
of $\rm{Spec}\ R$. More precisely, we know from \cite[Proposition
3.7]{BrGo} that for any Poisson prime, $P$, of $R$ the smooth
locus of the closed subvariety $\mathcal{V}(P)$ is a symplectic
leaf in $\rm{Spec}\ R$. This gives a one-to-one correspondence
between symplectic leaves in $\rm{Spec}\ R$ and Poisson prime
ideals of $R$: $\mathcal{S} \longleftrightarrow P$ where
$\overline{\mathcal{S}} = \mathcal{V}(P)$.

\subsection{}\label{proto}

Let $A$ be a $\C$-algebra. We shall say that a $\Z$-filtration,
$\F$, of $A$ is \textit{suitable} when we have:
$$ 0 = \F_{-1} \subseteq \F_{0} \subseteq \F_{1} \subseteq \dots
\subseteq A$$ and $\F$ satisfies (a) $\F_{i} \cdot \F_{j}
\subseteq \F_{i + j}$, (b) $\F_0 = \C$, (c) $\rm{dim}_{\C} \F_i <
\infty$ for all $i$, and (d) $\grA$ is an affine commutative
$\C$-algebra.

\begin{defn}

Let $A$ be a $\C$-algebra with suitable filtration, $\F$. We say
that $A$ has a \textit{proto-Poisson bracket with respect to $\F$}
if there exists a non-zero skew-symmetric $\C$-bilinear map
$\langle -,-\rangle : A \times A \to A$ which satisfies, for all
$a,b,c \in A$: \begin{enumerate}

\item $\langle ab,c \rangle = a\langle b,c
\rangle + \langle a,c \rangle b$.

\item $ \langle a, \langle b,c \rangle \rangle = \langle \langle
a,b \rangle, c\rangle + \langle b, \langle a,c \rangle \rangle$.

\item There is an integer $d$, the \textit{degree} of
$\langle -,-\rangle$, such that $\langle \F_i,\F_j \rangle
\subseteq \F_{i + j +d}$ for all $i,j \in \Z$, but there exist
$i,j \in \Z$ such that $\langle \F_i,\F_j \rangle \nsubseteq \F_{i
+ j + d - 1}$.
\end{enumerate}
\end{defn}

If the filtration is clear from the context we shall simply say
that $\langle -,- \rangle$ is a proto-Poisson bracket.

\begin{examples}{}
i) Let $A$ be an algebra with suitable filtration $\F$. Suppose
that $A$ is generated by $\F_1$ with $\F_i = (\F_{1})^i$, and that
$A$ is not commutative. Let $\langle -,- \rangle$ equal the
commutator bracket on $A$. Then $\langle -,- \rangle$ is a
proto-Poisson bracket on $A$. The only nontrivial condition to
check is (3). For this let $d$ be the integer such that $ \langle
\F_1, \F_1 \rangle \subseteq \F_{d+2}$ but $ \langle \F_1, \F_1
\rangle \nsubseteq \F_{d+1}$. There exists such an integer because
$A$ is not commutative. It is easily seen that $ \langle \F_i
,\F_j \rangle \subseteq \F_{i+j+d}$ for all $i,j \in \Z$.

ii) Let $R$ be a Poisson algebra with suitable filtration $\F$.
Let $\{-,-\}$ be the Poisson bracket on $R$; if $\{-,-\}$
satisfies condition (3) of the definition then it is a
proto-Poisson bracket. In particular, for any Poisson algebra,
$R$, with generating set $\{ a_1, \dots ,a_t\}$ we can define a
filtration, $\mathcal{F}$, by $\F_i = 0$ for $i < 0$, $\F_0 = \C$,
$\F_1 = \C 1 + \sum_{i = 1}^{t} \C a_i$ and $\F_i = ({\F_1})^i$
for $i \geq 2$. Then $\F$ is a suitable filtration and the Poisson
bracket on $R$ is a proto-Poisson bracket (in particular,
condition (3) will hold for some choice of $d$).
\\
\end{examples}

\begin{rem}
We note that for a commutative algebra $A$, the commutator bracket
is identically zero and so this is never an example of a
proto-Poisson bracket.
\end{rem}

We shall say that an ideal $I$ of $A$ is a
$\langle-,-\rangle$-\textit{ideal} if $\langle A,I\rangle
\subseteq I$. In example i), a  $\langle -,- \rangle$-ideal is
just an ideal of $A$; in example ii), a $\langle -,-\rangle$-ideal
is a Poisson ideal. The following extension to the present setting
of the standard Gabber-Hayashi recipe (see \cite[$\S 15$]{EG}, for
example) for constructing a Poisson bracket in Example (i) has a
routine proof which is left to the reader.

\begin{lem}

Let $A$ be a $\C$-algebra with suitable filtration, $\F$, and
proto-Poisson bracket, $\langle -,-\rangle$ of degree $d$. For
homogeneous elements $x,y \in \grA$ of degree $k$ and $l$
respectively, denote lifts of $x$ and $y$ by $\tilde{x}, \tilde{y}
\in A$, that is, $\sigma_{k}(\tilde{x}) = x$ and
$\sigma_{l}(\tilde{y}) = y$ where $\sigma$ denotes the principal
symbol map. Then
\begin{displaymath}
\mathrm{gr}\langle -,-\rangle : \grA \times \grA \to\grA;\ \ \
(x,y) \mapsto \sigma_{i + j + d}(\langle \tilde{x},
             \tilde{y}\rangle)
\end{displaymath}

\noindent defines a Poisson bracket of degree $d$ on $\grA$ when
extended
linearly.\\
\end{lem}

\subsection{} We can now state our main theorem.

\begin{thrm}\label{Ithrm}
Let $A$ be a $\C$-algebra with suitable filtration, $\F$, and
proto-Poisson bracket $\langle -,- \rangle$. Let $\grA$ have
Poisson bracket $\mathrm{gr}\langle -,-\rangle$ and let $I$ be a
prime $\langle -,-\rangle$-ideal of $A$. Suppose $X = \rm{Spec}\
\grA$ has FMSL with respect to the Poisson bracket induced on it
by $\langle -,- \rangle$. Let $V = \mathcal{V}(\grI)$. Then $V$ is
irreducible, and is the closure of a symplectic leaf in $X$.
\end{thrm}

The second claim follows quickly from the first. For suppose that
$V$ is irreducible. Since $I$ is a $\langle -,- \rangle$-ideal it
can easily be seen that $\grI$ is a Poisson ideal, and therefore
$\mathrm{rad}(\grI)$ is also Poisson by \cite[3.3.2]{D}. Hence $V$
is a closed irreducible Poisson subvariety and is the closure of a
symplectic leaf in $X$, as discussed in \ref{Poisson}. We note
that a version of Theorem \ref{Ithrm} is true with the weaker
assumption that $V$ (and not $X$) has FMSL. Then $V$ is
irreducible, but is not necessarily the closure of a symplectic
leaf of $X$.\\

\section{Applications}\label{appl}

\subsection{}We get as a corollary to Theorem \ref{Ithrm} a proof
of the result of Borho and Brylinski, and Joseph. There is a
detailed account, including the background material required, in
\cite[$\S 3-4$]{Vo}.

\begin{cor1}[\cite{BB}, \cite{J}]
Let $\mathfrak{g}$ be a complex semisimple Lie algebra, and let
$\mathcal{U} (\mathfrak{g})$ denote its enveloping algebra. Then
the associated variety of a primitive ideal of $\mathcal{U}
(\mathfrak{g})$ is irreducible.
\end{cor1}

\begin{proof} Let $\{ X_1, \dots ,X_m \}$ be a basis of
$\mathfrak{g}$ and let $[-,-]$ denote its Lie bracket. There is
suitable filtration, $\mathcal{B}$, on $\Ug$ where $\mathcal{B}_1
= \mathfrak{g}$ generates $\Ug$ as an algebra and $\mathcal{B}_i =
(\mathcal{B}_1)^i$ for all $i \geq 1$. Now $\mathrm{gr} \Ug =
\C[X_1, \dots ,X_m]$, and the variety $\mathrm{Spec}\ \mathrm{gr}
\Ug$ can be identified with $\mathfrak{g}^*$. As explained in
Examples \ref{proto} (i), setting $\langle -,- \rangle$ equal to
the commutator on $\Ug$ defines a proto-Poisson bracket on $\Ug$.
Therefore there is a Poisson bracket, $\rm{gr} \langle -,-
\rangle$, on $\mathrm{gr} \Ug$. However, since $\mathrm{gr} \Ug =
\C[X_1, \dots ,X_m]$, it is clear that $\rm{gr}\langle -,-
\rangle$ is extended from the Lie bracket on $\mathfrak{g}$,
giving the so-called Kostant-Kirillov Poisson bracket. Let $P$ be
a primitive ideal of $\mathcal{U}(\mathfrak{g})$, and let $Q$ be a
minimal primitive ideal contained in $P$. Now, $\mathfrak{g}^*$
will never have FMSL, but the Poisson subvariety $\mathcal{V}
(\mathrm{gr} Q)$ does (by \cite[Theorem 5.8]{Vo}) , the leaves
being the nilpotent coadjoint orbits (see \cite[Theorem
14.3.1]{Ma}). If we now take $A = \mathcal{U}(\mathfrak{g})/Q$
with filtration and proto-Poisson bracket induced from $\Ug$,
Theorem \ref{Ithrm} tells us that $\var(\mathrm{gr} P)$ is an
irreducible subvariety of $\var(\mathrm{gr} Q)$ and therefore also
of $\mathfrak{g}^*$.
\end{proof}

\subsection{}\label{v/G}

Before discussing further applications of Theorem \ref{Ithrm} we
introduce quotient varieties $\fv /G$ and describe their
symplectic leaves.\\

Let $(\fv, \omega, G)$ be an indecomposable symplectic triple (see
\cite[$\S 1$]{EG}) - in particular, $\fv$ is an even dimensional
$\C$-vector space, $\omega$ a symplectic form on $\fv$ and $G$ a
finite subgroup of the symplectic group of $\fv$. Let $\C[ \fv]$
be the coordinate ring of $\fv$ and let $\C[ \fv] * G$ be the skew
group ring. This latter algebra has centre $\C[ \fv]^G$, the ring
of $G$-invariants of $\C[ \fv]$. Let $e = \frac{1}{|G|} \sum_{g
\in G} g \in \C[G]$, then there is an isomorphism of
algebras $e \C[ \fv] * G e \cong \C[\fv]^G$.\\

The ring of invariants is a Poisson domain with bracket induced by
$\omega$ which we denote $\{-,-\}_{\omega}$. We note that
$\C[\fv]$ is a graded algebra and that $\{-,-\}_{\omega}$ has
degree $-2$. The variety $\fv /G = \Spec\ \C[\fv]^G$ has finitely
many symplectic leaves and, moreover, the leaves have been
described in \cite[Proposition 7.4]{BrGo}. Let $\pi : \fv \to \fv
/G$ be the orbit map and for $v \in \fv$ let $G_v$ denote the
stabiliser of $v$ in $G$. Given a subgroup $H$ of $G$ let $\fv_H^o
= \{ v \in \fv : H = G_v\}$. The symplectic leaves of $\fv /G$ are
the sets $\pi (\fv_H^o)$ as $H$ runs through subgroups of $G$ for
which $\fv_H^o \neq \emptyset$. If $H$ and $H'$ are conjugate
subgroups of $G$ then $\pi (\fv_H^o) = \pi (\fv_{H'}^o)$ so in
fact the leaves are in one-to-one correspondence with the
conjugacy classes of subgroups of $G$
which occur as the stabiliser of some element of $\fv$.\\

\subsection{Symplectic reflection algebras}\label{SRA}

For details of the following see \cite{EG}.

The symplectic reflection algebras corresponding to $(\fv, \omega,
G)$, written $H_{t, \mathbf{c}}$ where $t \in \C$ and $\mathbf{c}
\in \C^r$ for some $r$, are isomorphic, as vector spaces, to
$\C[\fv] \otimes_{\C} \C[G]$. They are deformations of the skew
group ring in the sense that, when they are filtered by putting
elements of $\fv$ in degree one, and putting $\C[G]$ in degree
zero, then the associated graded algebras are isomorphic to
$\C[\fv]*G$ (\cite[Theorem 1.3]{EG}). The spherical subalgebras $e
H_{t, \mathbf{c}} e$ inherit the filtration, and we denote this
filtration by $\mathcal{B}$. Their associated graded algebra is
$\C[\fv]^G$ (a consequence of $e$ being in degree zero). The
algebras $e H_{t,\mathbf{c}} e$ are commutative if and only if $t$
is zero (\cite[Theorem 1.6]{EG}).

\subsection{} We derive \cite[Theorem 2.1]{Gi} as a special case of Theorem
\ref{Ithrm}. We first require an elementary lemma.

\begin{lem}\label{2P}
Let $R$ be an affine commutative $\C$-algebra with two Poisson
brackets, $\{-,-\}_1$ and $\{-,-\}_2$, such that $\{-,-\}_1 =
\lambda \{-,-\}_2$ for some non-zero $\lambda \in \C$. Let $X =
\rm{Spec}\ R$. Then the symplectic leaves of $X$ with respect to
$\{-,-\}_1$ are the same as the symplectic leaves of $X$ with
respect to $\{-,-\}_2$. In particular, $(X, \{-,-\}_1)$ has FMSL
if and only if $(X,\{-,-\}_2)$ has FMSL.
\end{lem}

\begin{proof}
The rank of $\{-,-\}_1$ at any closed point $\mathbf{m}$ of $X$ is
equal to the rank of $\{-,-\}_2$ at $\mathbf{m}$, so the lemma
follows from the definition of symplectic leaf.
\end{proof}

\begin{cor2}
Let $eH_{t, \bf{c}}e$ be the algebra described in \ref{SRA}, and
let $t \neq 0$. Then for any primitive ideal $I$ of $eH_{t,
\bf{c}}e$, the variety $\mathcal{V} (\rm{gr}^{\mathcal{B}}I)$ is
irreducible.
\end{cor2}

\begin{proof}

 The filtration, $\mathcal{B}$, described in \ref{SRA} is a
suitable filtration on $eH_{t, \bf{c}}e$. Let $[-,-]$ be the
commutator bracket, then we claim that this is a proto-Poisson
bracket on $eH_{t, \bf{c}}e$. The only condition of definition
\ref{proto} which is non-trivial is (3), but this follows from
\cite[Claim 2.25(i)]{EG} (the degree, $d$, is $-2$ in this case).
Therefore $\rm{gr}[-,-]$ is a Poisson bracket of degree $-2$ and
so by \cite[Lemma 2.23(i)]{EG} there is some non-zero $\lambda \in
\C$ such that $\rm{gr}[-,-] = \lambda \{-,-\}_{\omega}$. By Lemma
\ref{2P}, $\rm{Spec}\ \C[\fv]^G$, with Poisson bracket
$\rm{gr}[-,-]$, has FMSL. We can now apply Theorem \ref{Ithrm}:
for any prime ideal $I$ of $eH_{t, \bf{c}}e$, $\mathcal{V}
(\rm{gr}^{\mathcal{B}}I)$ is irreducible. In particular, this is
true for any primitive ideal $I$.
\end{proof}

We extend this to the case when $t = 0$.

\begin{cor3}
Let $A = eH_{0,\bf{c}} e$ with filtration $\mathcal{B}$ as in
(\ref{SRA}) and denote its Poisson bracket by $\{-,-\}$ . Let $I$
be a Poisson prime ideal of $A$, then
$\mathcal{V}(\rm{gr}^\mathcal{B} I)$ is irreducible.
\end{cor3}

\begin{proof}
The filtration $\mathcal{B}$ is suitable. We would like $\{-,-\}$
to be a proto-Poisson bracket, and so, as noted in Examples
\ref{proto} ii), we need to show that condition (3) of Definition
\ref{proto} is satisfied. It can be seen from the construction of
$\{-,-\}$ that there is some $l \geq 2$ so that $
 \{ \F_{i}, \F_{j}\} \subseteq \F_{i + j - l},\
\rm{for\ all}\ i,j \in \Z$ but $ \{ \F_{i}, \F_{j}\} \nsubseteq
\F_{i + j - (l + 1)},\ \rm{for\ some}\ i,j \in \Z.$ By \cite[Lemma
2.26]{EG}, we must have $l = 2$.

It remains to show that $\Spec\ \rm{gr}^{\mathcal{B}} A = \fv/G$
with Poisson bracket $\rm{gr} \{-,-\}$ has FMSL. However, by
\cite[Lemma 2.23]{EG}, $\rm{gr} \{-,-\} = \lambda
\{-,-\}_{\omega}$ for some $\lambda \in \C^{*}$. Therefore by
\ref{v/G} and Lemma \ref{2P} , $(\fv/G, \rm{gr} \{-,-\})$ has
FMSL. Hence $\var (\rm{gr}^{\mathcal{B}}I)$ is irreducible.
\end{proof}

\subsection{}\label{rep}

Let $Z = Z_{0 , \mathbf{c}}$ and let $H = H_{0,\bf{c}}$. Our
objective in proving Theorem \ref{Ithrm} is to better understand
the symplectic leaves of $\rm{Spec}\ Z$. For it was shown in
\cite[Theorem 4.2]{BrGo} that the symplectic leaves of
$\mathrm{Spec} Z$ control the finite dimensional representation
theory of the corresponding symplectic reflection algebra. In more
detail, \cite[Theorem 4.2]{BrGo} says that if two closed points
$\bf{m}$ and $\bf{n}$ of $\rm{Spec}\ Z$ lie in the same symplectic
leaf then $H/{\mathbf{m} H}$ and $H/{\mathbf{n} H}$ are isomorphic
$\C$-algebras.

By \cite[Theorem 7.8]{BrGo}, $\rm{Spec}\ Z$ has FMSL. Therefore we
can restate our goal as finding a description of the Poisson prime
ideals of $Z$. By \ref{v/G}, the corollary below allows us to
attach to each Poisson prime of $Z$, a conjugacy class of a
certain subgroup of $G$.

\subsection{}\label{t=0} \textbf{The case $\bf{t = 0}$:}  Let $Z$
and $H$ be as in \ref{rep}. The algebras $e H_{0, \mathbf{c}} e$
and $Z$ are both Poisson algebras via the Gabber-Hayashi
construction.

\begin{thrm}[\cite{EG},Theorem 3.1]
The map
\begin{displaymath}
\psi: Z \to e H_{0, \mathbf{c}} e
\end{displaymath}
\begin{displaymath}
z \mapsto eze
\end{displaymath}
is a Poisson isomorphism.
\end{thrm}

Let $\mathcal{A}$ denote the filtration on $Z$ induced from that
on $H$. The map $\psi$ preserves the filtrations $\mathcal{A}$ and
$\mathcal{B}$ and we have the following result.

\begin{prop}[\cite{EG}, Proposition 3.4]
The associated graded map $\rm{gr(\psi)} : \rm{gr}^{\mathcal{A}}Z
\to \rm{gr}^{\mathcal{B}}e H_{0, \mathbf{c}} e$ is an algebra
isomorphism.
\end{prop}

Let $P$ be a prime ideal of $Z$, then by the Theorem, $P$ is
Poisson if and only if $\psi(P)$ is Poisson. Furthermore, by the
Proposition, $\var (\rm{gr}^{\mathcal{A}} P)$ is irreducible if
and only if $\var (\rm{gr}^{\mathcal{B}} \psi(P))$ is irreducible.

\begin{cor4}
Let $I$ be a Poisson prime ideal of $Z$. Then $\var
(\rm{gr}^{\mathcal{A}} I)$ is irreducible.
\end{cor4}
\begin{proof}
This is immediate from the previous paragraph and Corollary 3.
\end{proof}

\section{Preliminaries to the proof of theorem \ref{Ithrm}}\label{Prelim}

\subsection{} For the remainder of the paper retain the following
notation: Let $A$, $\mathcal{F}$, $\langle -,- \rangle$, $I$, $X$
and $V$ satisfy all of the hypotheses of Theorem \ref{Ithrm} and
let $M = A/I$. We can choose an irreducible component of $V$ of
maximal dimension. As explained in \ref{Poisson}, because $X$ has
FMSL, there exists a symplectic leaf $\s$ such that $\os$ is this
component. Then dim$\s$ $=$ dim$\os = n$, and by definition, the
closed points of $\s$ all have rank $n$. For a subvariety $W$ of
$V$ we write $\mathrm{sm}W$ for the smooth locus of $W$ and
$\mathrm{sing}W$ for $W \setminus \mathrm{sm}W$.

\begin{lem}\label{SLlem}

$\mathrm{(1)}$ $\s$ is open in $V$.\\

$\mathrm{(2)}$ \textrm{dim} $ \os\setminus\s \leq \rm{dim}\os - 2 = \rm{dim}V - 2$.\\

$\mathrm{(3)}$ $\s$ is a homogeneous subvariety of $V$ i.e. there
exist homogeneous elements\\ \noindent $ g_{1}, \dots ,g_{s} \in
{\grA}/{\grI}$ such that $\s = V\setminus \var(g_{1}, \dots
,g_{s})$.

\end{lem}

\begin{proof}

$\mathrm{(1)}$ Let $V = I_1\cup \dots \cup I_k$ be an irredundant
irreducible decomposition of $V$ with $I_1 = \os$. We claim that
$\s \cap I_j = \emptyset$ for all $2\leq j\leq k$. If not, then
for some $j$, $\s \cap I_j$ contains a closed point of rank $n$,
$\mathbf m$ say. By \cite[Proposition 3.7]{BrGo} the smooth locus
of $I_j$, $\mathrm{sm}I_j$, is a symplectic leaf in $V$ which
contains $\mathbf m$. Then $\s = \mathrm{sm}I_j$ implies that $I_1
= I_j$, a contradiction.\smallskip

Therefore $V\setminus \s = (\os \setminus \s)\cup I_2\cup \dots
\cup I_k$ is closed in $V$.\\

$\mathrm{(2)}$ It is clear from the definition of symplectic
leaves that they are even dimensional. We can write $\os \setminus
\s$ as a finite union of symplectic leaves, each of which has
dimension less than $\s$. The
inequality follows because of even dimensionality.\\

$\mathrm{(3)}$ This is true because $\s = \mathrm{sm}I_1$. Since
$V$ is homogeneous $I_1$ is also homogeneous so we may assume that
the ideal of $I_1$ is generated by homogeneous elements $h_1,
\dots , h_l$ in some polynomial ring $\C [x_1, \dots x_m]$. Now
$\mathrm{sing}I_1$ is defined as the points vanishing at certain
$(m - r) \times (m-r)$ minors of
\begin{displaymath} \left(
\begin{array}{c} {\partial h_i}/{\partial
x_j}
\end{array} \right).
\end{displaymath}
These minors are homogeneous polynomials in the $x_j
\textstyle{s}$. Hence $\mathrm{sing}I_1$ and also $\mathrm{sm}I_1$
are homogeneous subvarieties.

\end{proof}

\noindent (1) and (3) imply that there is a homogeneous open set
$U \subseteq X$ such that $U \cap V = \s$. We write $U =
X\setminus \var(f_1, \dots ,f_t )$ where the $f_i$ are homogeneous
elements of $\grA$.

\subsection{Microlocalisation}\label{mic}

The proof of the main theorem makes use of microlocalisation
techniques which are described in \cite[$\S$ 3-4]{Vo}. We shall
say that a $\Z$-filtration, $\mathcal{B}$, of an $A$-module $M$ is
\textit{good} if $\bigcap_{n \in \Z} \mathcal{B}_n = 0$,
$\bigcup_{n \in \Z} \mathcal{B}_n = M$ and
$\mathrm{gr}^{\mathcal{B}} M$ is a finitely generated $\grA
$-module. We recall the definition of support. Let $R$ be a
commutative ring and let $N$ be an $R$-module. Then
\begin{displaymath}
\mathrm{supp}_{R} N = \{ P \in \mathrm{Spec}\ R : N_{P} \neq 0 \}.
\end{displaymath}

\noindent We write $\rm{supp}\ N$ when the ring is clear from the
context.

For $M = A/I$, consider the induced filtration on $M$ (which we
also call $\mathcal{F}$). Then $$ \mathrm{supp}\ \grM = \var
(\mathrm{ann}\ \grM) = \var (\grI ) = V$$ where the first equality
is true because $\mathcal{F}$ is a good filtration of $M$. Now
$\grM$ defines a sheaf of $\mathcal{O}_X$-modules, $\mathcal{M}$
on $X$. We only need to know the sections of $\mathcal{M}$ over
$U$, which we can calculate explicitly:

\begin{displaymath}
\mathcal{M}(U) = \{ (m_{f_i}) \in \prod_{i=1}^{t}\ (\grM )_{f_i} :
m_{f_i} = m_{f_j} \in (\grM )_{f_{i}f_{j}} \forall i,j \} .
\end{displaymath}

\begin{lem}[\cite{Vo}, Lemma 3.3]\label{V1}
Let $K$ be the kernel of the natural $\grA$-module map
\begin{displaymath}
\beta : \grM \rightarrow \mathcal{M}(U)
\end{displaymath}
\begin{displaymath}
m \mapsto (m).
\end{displaymath}

\noindent Then\ $K = \{m\in \grM : \mathrm{for\ each}\ i\
\mathrm{there\ exists}\ N_i \in \N\ \mathrm{such\ that}\
{f_{i}}^{N_{i}}m = 0 \}$, and
\newline $ \mathrm{supp} \ K \cap U = \emptyset $.

\end{lem}

Microlocalisation introduces a new filtration, $\Gamma$, on $M$
which is compatible with $\mathcal{F}$ \cite[Corollary 6.9]{Vo}.
We can give a description of $\Gamma$ in terms of the $f_i
\textstyle{s}$ introduced above. Let $p_i$ be the degree of $f_i$
and for each $i$ choose a lift, $\phi_i$, of $f_i$ to $A$. Thus
$\phi_i \in \F_{p_i}$ and $\sigma_{p_i} (\phi_i) = f_i$.

Let $\mathcal{I} = \{1, \dots , t\}$ and suppose $\tau = (i_1,
\dots ,i_N) \in \mathcal{I}^N$ is an ordered $N$-tuple of elements
of $\mathcal{I}$. Define

\begin{displaymath}
       p_{\tau} = \sum_{j=1}^N p_{i_j},\ \ \  \phi_{\tau} =
       \prod_{j=1}^N \phi_{i_j} \in \F_{p_{\tau}}.
\end{displaymath}

\noindent Then $ \Gamma_n = \{ m \in M :\ \mathrm{for\ all}\ N\
\mathrm{sufficiently\ large,\ and\ for\ all}\ \tau \in
\mathcal{I}^N, \phi_{\tau} \cdot m \in \F_{n+p_{\tau}} \}.$ We see
$\Gamma$ has the property that $\mathcal{F}_i \subseteq \Gamma_i\
\forall\ i \in \Z$ so that $\grGM$ is a $\grA$-module and there is
a canonical map (of $\grA$-modules) $\alpha : \grM \rightarrow
\grGM$.

\begin{prop}[\cite{Vo}, Proposition 3.11]\label{commD}
There\ exists\ a\ map\ of\ $\grA$-modules $\theta: \grGM \to
\mathcal{M}(U)$ which is injective and gives rise to the following
commutative diagram:

$$\diagram \grM \drto_{\beta} \rto^{\alpha} & \grGM \dto^{\theta} \\
                                            & \mathcal{M}(U)
                                            \enddiagram $$
\end{prop}

\subsection{}

 For a commutative ring $R$ and an $R$-module $N$ we
recall that the \textit{associated primes} of $N$, written Ass$N$,
are the set of primes of $R$ which are annihilators of elements of
$N$. The following result will be key to the proof of Theorem
\ref{Ithrm}: we will use it to prove that $\mathcal{M} (U)$ is a
finitely generated $\grA$-module. In fact, we show later that if
we take $R = \grA$, $N = \grM$ and $W = U$ in the statement below
then condition (\ref{Eq2}) is a consequence of Lemma \ref{SLlem}
(2).

\begin{thrm}[\cite{Gr}, Proposition 5.11.1]\label{Gr}
Suppose $R$ is an affine commutative algebra over $\C$, $N$ is a
finitely generated $R$-module and $W$ is an open set in
$\mathrm{Spec} R$. Let $\mathcal{N}$ denote the sheaf of modules
associated to $N$. Then the $R$-module $\mathcal{N}(W)$ is
finitely generated if and only if for every prime $P \in W \cap
\mathrm{Ass} N$, $\overline{P}$, the closure of $P$ in
$\mathrm{Spec} R$, satisfies

\numberwithin{equation}{section}
\begin{equation}\label{Eq2}
\overline{P} \cap (\mathrm{Spec}R \backslash W)\ \mathrm{has\
codimension\ at\ least\ 2\ in}\ \overline{P}.
\end{equation}
\end{thrm}

We shall also require some information about the associated primes
of $\grM$. Let $R$ be a Poisson algebra and $N$ an $R$-module,
then we say that $N$ is a \textit{Poisson module} if there exists
a $\C$-bilinear form $\{-,-\}_N : R \times N \to N$ satisfying
$\{r,r'n\}_N = \{r,r'\}n + r'\{r,n\}_N$ for all $r, r' \in R$ and
$n \in N$. It is clear that, in our setting, $\grM$ is a Poisson
$\grA$-module.

\begin{lem}[\cite{CO}, Theorem 4.5]\label{ass}
Let $R$ be a Noetherian Poisson algebra and $N$ be a finitely
generated Poisson $R$-module. Then the associated primes of $N$
are Poisson ideals of $R$.
\end{lem}

\begin{proof}
Let $P$ be an associated prime of $N$. Let $L = \{ n \in N: P^i n=
0\ \mathrm{for\ some}\ i \geq 0\}$, this is a non-zero submodule
of $N$. We claim that $L$ is a Poisson submodule of $N$. To show
this we first note that since $L$ is finitely generated there is
some $t \geq 1$ such that $P^t L = 0$. Let $l \in L, r \in R$. For
all $r' \in P^t$, $$ 0 = \{r,r'l\}_N = \{r,r'\}l + r'\{r,l\}_N.$$
Therefore $ P^t\{R,L\}_N \subseteq \{R,P^t\}L \subseteq L,$ which
implies that $P^{2t}\{R,L\}_N = 0$. Hence $\{R,L\}_N \subseteq L$,
by definition of $L$, which means that $L$ is a Poisson submodule
of $N$. By \cite[Lemma 4.1]{BrGo}, $\mathfrak{I} = \mathrm{ann}_R
L$ is a Poisson ideal of $R$. There is some element $x \in L$ such
that $\mathrm{Ann}_R\{x\} = P$, then $x \in L$ implies
$\mathfrak{I} \subseteq P$. Taking radicals of the ideals $P^t
\subseteq \mathfrak{I} \subseteq P$ yields
$\mathrm{rad}\mathfrak{I} = P$, and therefore $P$ is a Poisson
ideal by \cite[3.3.2]{D}.
\end{proof}

\bigskip

\section{Proof of Theorem \ref{Ithrm}}\label{proof}

\subsection{} We retain the notation introduced at the beginning
of section \ref{Prelim}, in particular we recall that $M = A/I$.
We also use the notation $\alpha, \beta$ and $\theta$ from
Proposition \ref{commD}.

We make the following two assumptions:

\begin{claim}[\bf{1}]
$\Gamma$ is a good filtration of $M$.
\end{claim}

\begin{claim}[\bf{2}]
$\mathrm{supp}_{\grA}\ \mathcal{M}(U) \subseteq \os$.
\end{claim}

Now $\mathrm{supp}_{\grA}\ \grGM \subseteq \mathrm{supp}_{\grA}\
\mathcal{M}(U)$ because $\theta$ is injective. The left hand side
equals $V$ by claim (1) and the right hand side is contained in
$\os$ by claim (2). So we have $V \subseteq \os$ which implies
that $V = \os$ and this proves the theorem.
\\

It remains to prove the two claims.\\

\noindent {\it{Proof of Claim (1)}}: Recall that
there are three conditions to check.\\

(a) $\mathit{\bigcap_{n \in \Z} \Gamma_{n} = 0}$. Let $M_{-
\infty} = \bigcap_{n \in \Z} \Gamma_{n}$. It is easy to check that
$M_{- \infty}$ is an $A$-sub-bimodule of $M$. We see that
$\mathrm{gr}^{\Gamma} (M_{- \infty}) = 0$: for all $i \in \Z$,
$(\Gamma_i \cap M_{- \infty})/(\Gamma_{i-1} \cap M_{- \infty}) =
(\Gamma_{i-1} \cap M_{- \infty})/(\Gamma_{i-1} \cap M_{- \infty})
= 0$. Therefore the map $\mathrm{gr}^{\mathcal{F}} (M_{- \infty})
\to \mathrm{gr}^{\Gamma} (M_{- \infty})$ given by the restriction
of the map $\alpha$ above, is the zero map. It follows from
Proposition \ref{commD} that $\mathrm{gr}^{\mathcal{F}} (M_{-
\infty}) \subseteq K$ where $K$ is the kernel of $\beta$. By Lemma
\ref{V1}
\begin{equation}\label{Grlem}
 \mathrm{supp}\ \mathrm{gr}^{\mathcal{F}} (M_{- \infty}) \cap U =
 \emptyset.
\end{equation}

\par \noindent Now since $M_{- \infty}$ is an $A$-sub-bimodule of $A/I$, there is an ideal
$J$ of $A$ such that $I\subseteq J\subseteq A$ and $M_{- \infty} =
J/I$. Suppose that $M_{- \infty} \neq 0$. Then $J$ properly
contains $I$. It is a consequence of \cite[Propositions 3.15 and
6.6]{KL} that $\mathrm{dim} \var (\mathrm{gr}^{\mathcal{F}} J)$ $<
\mathrm{dim} \var (\mathrm{gr}^{\mathcal{F}} I) = \mathrm{dim}\os
$. Let $\mathfrak{p}$ be the defining ideal of $\os$. The equality
of closed sets (where the support is considered over $\grA$)
$$ \mathrm{supp}\ \grf (A/I) = \mathrm{supp}\ \grf (J/I) \cup
\mathrm{supp}\ \grf (A/J)$$ implies that $\mathfrak{p} \in
\mathrm{supp}\ \grf (J/I)$ and therefore that $\os \subseteq
\mathrm{supp}\ \grf (J/I)$. Hence $\s \subseteq \os \subseteq
\mathrm{supp}\ \grf (J/I) = \mathrm{supp}\ \grf (M_{- \infty})$.
This contradicts (\ref{Grlem}) and so $M_{- \infty} =
0$.\\

(b) $\mathit{ \bigcup_{n\in \Z} \Gamma_n = M}$. This is
straightforward because $\mathcal{F}_n \subseteq \Gamma_n$ implies
$M = \bigcup_{n} \mathcal{F}_n \subseteq \bigcup_{n} \Gamma_n
\subseteq M$.\\

(c) $\mathit{\grGM\ is\ a\ finitely\ generated\ \grA-module}$. To
prove this we in fact show that $\mathcal{M}(U)$ is a finitely
generated $\grA$-module (which proves (c) by Proposition
\ref{commD}, since $\grA$ is Noetherian). We would like to show
that $\mathcal{M}(U)$ is finitely generated and so by Theorem
\ref{Gr} it suffices to show that each prime $P \in U \cap
\mathrm{Ass}\ \grM$ satisfies (\ref{Eq2}) with $R = \grA$. Let $P
\in U \cap \mathrm{Ass}\ \grM \subseteq U \cap \mathrm{supp}\ \grM
= \s $. By Lemma \ref{ass}, $P$ is a Poisson prime ideal of
$\grA$. Now $U \cap \overline{P}$ is a nonempty open subset of
$\overline{P}$ which means that it contains a closed point
$\mathbf{m}$, and $U \cap \overline{P} \subseteq \s$ implies that
$\mathbf{m}$ has rank $n$. We conclude that $\mathrm{dim}
\overline{P} = n$ by \cite[Lemma 3.1(5)]{BrGo}, and therefore
$\overline{P} = \os$. Thus condition (\ref{Eq2}) is a consequence
of Lemma \ref{SLlem} (2) and $\mathcal{M}(U)$ is finitely
generated. This proves (c) and concludes the proof of
Claim (1).\\

\noindent $\it{Proof\ of\ Claim\ (2)}$: We must show that $P
\notin \os \Rightarrow \mathcal{M}(U)_{P} = 0$. Let $P \notin
\os$. This is equivalent to
 there being a neighbourhood, $Y$, of $P$ in $X$ such that $Y \cap \s
= \emptyset$. Without loss of generality, we may assume that $Y$
is some standard open set $O_{g}$ for some $g \in \grA$ with $g
\notin P$. We have $U = X \setminus \mathcal{V}(f_1, \dots ,f_t) =
O_{f_1} \cup \dots \cup O_{f_t}$; for each $i$, $O_{f_{i} g} =
O_{f_{i}} \cap O_g$ is an open subset of $U$ which intersects $\s$
trivially. Since $U \cap V = \s$ we conclude that $O_{f_{i}g}$ is
contained in the open set $X \setminus V$, and therefore that $V
\subseteq \mathcal{V}(f_{i}g)$. By considering the ideals of these
subvarieties we deduce that $\mathrm{rad} \langle f_{i} g\rangle
\subseteq \mathrm{rad} (\mathrm{ann}\ \grM)$. Hence there are
integers $k_{i}$ so that $(f_{i} g)^{k_i} \in \mathrm{ann}\ \grM$.
We consider a typical element $(\frac{m_1}{{f_1}^{N_1}}, \dots
,\frac{m_t}{{f_t}^{N_t}}) \in \mathcal{M}(U)$. Let $k = max_i \{
k_i \}$, then for all $i$:
\begin{displaymath} (f_i g)^k m_i = 0\\
\Rightarrow \frac{g^k m_i}{{f_i}^{N_i}} = 0 \in \grM_{f_i}\\
\Rightarrow \mathcal{M}(U)_P = 0
\end{displaymath}
and this proves Claim (2).

\bibliographystyle{amsplain}

\begin{thebibliography}{10}

\bibitem{BB}
W.~Borho and J-L. Brylinski, \emph{Differential operators on
homogeneous spaces
  iii}, Inventiones Math. \textbf{80} (1985), 1--68.

\bibitem{BrGo}
K.A. Brown and I.~Gordon, \emph{Poisson orders, symplectic
reflection algebras
  and representation theory}, J. reine angew. Math. \textbf{559} (2003),
  193--216.

\bibitem{CO}
C.~Casselman and M.~Osborne, \emph{The restriction of admissable
  representations to $\mathfrak{n}$}, Math. Ann. \textbf{233} (1978), 193--198.

\bibitem{D}
J.~Dixmier, \emph{Enveloping algebras}, Grad. Stud. Math.,
Vol.~11, Amer. Math.
  Soc. Providence, RI, 1996.

\bibitem{EG}
P.~Etingof and V.~Ginzburg, \emph{Symplectic reflection algebras,
  Calogero-Moser space, and deformed Harish-Chandra homomorphism}, Invent. Math
  \textbf{147} (2002), no.~2, 243--348.

\bibitem{Gi}
V.~Ginzburg, \emph{On primitive ideals}, Selecta Math. \textbf{New
ser. 9}
  (2003), 379--407.

\bibitem{Gr}
A.~Grothendieck and J.~Dieudonn\'{e}, \emph{El\'{e}ments de
g\'{e}om\'{e}trie alg\'{e}brique iv.
  etude locale des sch\'{e}mas et des morphismes de sch\'{e}mas (seconde partie)},
  Publications Math\'{e}matiques \textbf{24}, Institut des Hautes Etudes Scientifiques (1965).

\bibitem{J}
A.~Joseph, \emph{On the associated variety of a primitive ideal},
J. Algebra
  \textbf{93} (1985), 509--523.

\bibitem{KL}
G.R. Krause and T.H. Lenagan, \emph{Growth of algebras and
gelfand-kirillov
  dimension}, revised ed., Grad. Stud. Math., vol.~22, Amer. Math. Soc.
  Providence, RI, 2000.

\bibitem{Ma}
J.E. Marsden and T.S. Ratiu, \emph{Introduction to mechanics and
symmetry},
  second ed., Texts in Applied Mathematics, vol.~17, Springer-Verlag, New York,
  1999.

\bibitem{Vo}
D.A. Vogan, \emph{Associated varieties and unipotent
representations}, Progr.
  Math. \textbf{101} (1991), 315--388.

\bibitem{W}
A. Weinstein, \emph{The local structure of Poisson manifolds}, J.
Diff. Geom. \textbf{18} (1983),
 523--557.

\end{thebibliography}

\end{document}